\documentclass[12pt]{article}
\usepackage{amsfonts}
\usepackage{a4wide}
\usepackage{amsmath,amscd,amsthm,a4,amssymb}
\usepackage{graphicx}

\begin{document}

\begin{center}
{\Large \textbf{A note on the problem of straight-line interpolation by
ridge functions}}

\

\textbf{Azer Akhmedov}$^{1}$ \textbf{and Vugar E. Ismailov}\footnote{%
Corresponding author}$^{2,3}$ \vspace{3mm}

$^{1}${Department of Mathematics, North Dakota State University, Fargo, ND,
58102, USA}\vspace{1mm}

$^{2}${Institute of Mathematics and Mechanics, Baku, Azerbaijan} \vspace{1mm}

$^{3}${Center for Mathematics and its Applications, Khazar University, Baku,
Azerbaijan} \vspace{1mm}

e-mails: $^{1}$azer.akhmedov@ndsu.edu, $^{2}$vugaris@mail.ru
\end{center}

\smallskip

\textbf{Abstract.} This paper proves that any three distinct straight lines
in the plane contain vertices of a closed broken line with sides parallel to
the coordinate axes. We use this geometric fact to analyze the interpolation
problem on lines by linear combinations of ridge functions. We
constructively prove that it is impossible to interpolate arbitrary data on
any three or more straight lines with sums of ridge functions having only
two fixed directions. The general case, involving more directions and
arbitrarily many straight lines, is reduced to the problem of existence of
certain finite sets in the union of these lines.

\bigskip

\textit{Mathematics Subject Classifications:} 15A06, 26B40, 41A05

\textit{Keywords:} ridge function; interpolation; linear equation; closed
path; cycle

\bigskip

\begin{center}
{\large \textbf{1. Introduction}}
\end{center}

The current paper proves that any three distinct straight lines in the plane
always contain vertices of a closed broken line with sides parallel to the
coordinate axes. We show that this geometric fact is applicable to the ridge
function interpolation problem.

A \textit{ridge function} is a multivariate function of the format%
\begin{equation*}
F(\mathbf{x})=g\left( \mathbf{a}\cdot \mathbf{x}\right) =g\left(
a_{1}x_{1}+\cdot \cdot \cdot +a_{n}x_{n}\right) ,
\end{equation*}%
where $g:\mathbb{R}\rightarrow \mathbb{R}$ and $\mathbf{a}=\left(
a_{1},...,a_{n}\right) $ is a fixed vector (direction) in $\mathbb{R}%
^{n}\backslash \left\{ \mathbf{0}\right\} .$ In other words, a ridge
function is a multivariate function constant on the hyperplanes $\mathbf{a}%
\cdot \mathbf{x}=c$, for any $c\in $ $\mathbb{R}$. These functions arise
naturally in various fields. They arise in partial differential equations
(where they are called \textit{plane waves}, see, e.g., \cite{John}), in
computerized tomography (see, e.g., \cite{Log}), in statistics (especially,
in the theory of projection pursuit and projection regression; see, e.g.,
\cite{Fri}), in the theory of neural networks (see, e.g., \cite{Ism}).
Finally, these functions are used in modern approximation theory as an
effective and convenient tool for approximating and interpolating
complicated multivariate functions (see, e.g., \cite{BP,Mai}). For more on
ridge functions and application areas see the monographs by Pinkus \cite{Pin}
and Ismailov \cite{Ism}.

Let $\mathbf{a}^{1},...,\mathbf{a}^{r}$ be fixed pairwise linearly
independent directions in $\mathbb{R}^{n}\backslash \left\{ \mathbf{0}%
\right\} $. Consider the following set of linear combinations of ridge
functions

\begin{equation*}
\mathcal{R}\left( \mathbf{a}^{1},...,\mathbf{a}^{r}\right) =\left\{
\sum\limits_{i=1}^{r}g_{i}\left( \mathbf{a}^{i}\cdot \mathbf{x}\right)
:g_{i}:\mathbb{R}\rightarrow \mathbb{R},i=1,...,r\right\} .
\end{equation*}%
Note that the set $\mathcal{R}\left( \mathbf{a}^{1},...,\mathbf{a}%
^{r}\right) $ is a linear space. We are interested in the problem of
interpolation by functions from $\mathcal{R}\left( \mathbf{a}^{1},...,%
\mathbf{a}^{r}\right) $. Interpolation at a finite number of points by such
functions has been considered in Braess, Pinkus \cite{BP}, Sun \cite{Sun},
Reid, Sun \cite{Rei}, Weinmann \cite{Wei} and Levesley, Sun \cite{Lev}. The
results therein are complete only in the case of two directions ($r=2$), and
three directions in $\mathbb{R}^{2}$.

In \cite{IP}, Pinkus and Ismailov discussed the problem of interpolation by
functions from $\mathcal{R}\left( \mathbf{a}^{1},...,\mathbf{a}^{r}\right) $
in straight lines. To make the problem precise, assume we are given the
straight lines $\{t\mathbf{b}^{i}+\mathbf{c}^{i}:t\in \mathbb{R}\}$, $%
\mathbf{b}^{i}\neq 0,$ $i=1,...,m$. The question they asked is when, for
every choice of data $g_{i}(t)$, $i=1,...,m$, there exists a function $G\in
\mathcal{R}\left( \mathbf{a}^{1},...,\mathbf{a}^{r}\right) $ satisfying

\begin{equation*}
G(t\mathbf{b}^{i}+\mathbf{c}^{i})=g_{i}(t),
\end{equation*}%
for all $t\in \mathbb{R}$ and any $i=1,...,m$. The paper \cite{IP} totally
analyzed this question with one and two directions ($r\leq 2$). It was found
necessary and sufficient conditions for the possibility of interpolation on
one and two straight lines. One of the major results of \cite{IP} dealt with
the case of three (or more) distinct straight lines.

\bigskip

\textbf{Theorem 1.1.} (see \cite{IP}) \textit{Assume we are given linearly
independent directions $\mathbf{a}^{1},\mathbf{a}^{2}$ in $\mathbb{R}^{n}$.
Then for any three distinct straight lines $\{t\mathbf{b}^{i}+\mathbf{c}%
^{i}:t\in \mathbb{R}\}$, $i=1,2,3,$ and almost all $g_{1},g_{2},g_{3}$ defined on $%
\mathbb{R}$ there does not exist a $G\in \mathcal{R}\left( \mathbf{a}^{1},%
\mathbf{a}^{2}\right) $ satisfying}

\begin{equation*}
G(t\mathbf{b}^{i}+\mathbf{c}^{i})=g_{i}(t),~t\in \mathbb{R},~i=1,2,3.
\end{equation*}

\bigskip

Theorem 1.1 was proved by analyzing certain first order difference
equations. In this paper we give a completely different, geometrically
explicit and constructive proof for Theorem 1.1.

Our approach will be based on so called paths. A \textit{path} with respect
to the directions $\mathbf{a}^{1}$ and $\mathbf{a}^{2}$ is an ordered set of
points $(\mathbf{x}^{1},...,\mathbf{x}^{k})$ in $\mathbb{R}^{n}$ with $%
\mathbf{x}^{i}\neq \mathbf{x}^{i+1}$ and the units $\mathbf{x}^{i+1}-\mathbf{%
x}^{i}$ perpendicular alternatively to $\mathbf{a}^{1}$ and $\mathbf{a}^{2}$%
. That is, after a suitable permutation we have $\mathbf{a}^{1}\cdot $ $%
\mathbf{x}^{1}=\mathbf{a}^{1}\cdot $ $\mathbf{x}^{2},$ $\mathbf{a}^{2}\cdot $
$\mathbf{x}^{2}=\mathbf{a}^{2}\cdot $ $\mathbf{x}^{3},$ $\mathbf{a}^{1}\cdot
$ $\mathbf{x}^{3}=\mathbf{a}^{1}\cdot $ $\mathbf{x}^{4}$, and so on. If $%
\left( \mathbf{x}^{1},...,\mathbf{x}^{k},\mathbf{x}^{1}\right) $ is a path
and $k$ is an even number, then the path $\left( \mathbf{x}^{1},...,\mathbf{x%
}^{k}\right) $ is said to be a \textit{closed path}. For example, a set of
vertices (ordered clockwise or counterclockwise) of any polygon in $\mathbb{R%
}^{2}$ with sides parallel to the coordinate axis is a closed path with
respect to the coordinate directions $(1,0)$ and $(0,1)$. Fig. 1 describes
two different paths in the coordinate plane.

\begin{figure}[ht]
\hfill
\par
\begin{center}
\includegraphics[width=0.95\textwidth]{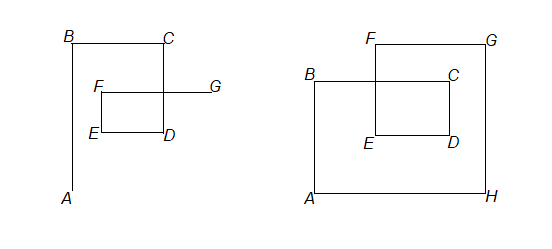}
\end{center}
\caption{\textit{In the first picture, the points $A,B,C,D,E,F,G$, in the
given order, form a path with respect to the coordinate directions. The
second picture illustrates a closed path.}}
\end{figure}

Paths are geometrically explicit objects. In the special case when $\mathbf{a%
}_{1}$ and $\mathbf{a}_{2}$ coincide with the coordinate directions these
objects were first introduced by Diliberto and Straus \cite{Dil} (in \cite%
{Dil}, they are called \textquotedblleft permissible lines\textquotedblright
). They appeared further in a number of papers with several different names
such as \textquotedblleft bolts" (see, e.g., \cite{Arn,Kh}),
\textquotedblleft trips\textquotedblright\ (see \cite{Mar2}),
\textquotedblleft links" (see, e.g., \cite{Cow}). In \cite{Is1}, paths were
generalized to those with respect to a finite set of functions. The last
objects turned out to be useful in problems of representation by linear
superpositions. For a detailed history of paths, their generalizations and
various applications see \cite{Ism}.

It should be remarked that the problem of interpolation of arbitrary data on
the lines $l_{i}=\{t\mathbf{b}^{i}+\mathbf{c}^{i}:t\in \mathbb{R}\}$, $%
i=1,...,m,$ by functions from $\mathcal{R}\left( \mathbf{a}^{1},...,\mathbf{a%
}^{r}\right) $ is equivalent to the problem of the representation of an
arbitrarily given $F$ defined on the union $L=\cup _{i=1}^{m}l_{i}$ by such
functions (see \cite{IP}). Concerning this problem, we have the following
result from \cite{Is1}.

\bigskip

\textbf{Theorem 1.2.} \textit{Assume we are given two directions $\mathbf{a}%
^{1}$ and $\mathbf{a}^{2}$ in $\mathbb{R}^{n}$. Let $X$ be any subset of $%
\mathbb{R}^{n}$. Then every given function defined on $X$ is in $\mathcal{R}%
\left( \mathbf{a}^{1},\mathbf{a}^{2}\right) $ if and only if there are no
finite set of points in $X$ that form a closed path with respect to the
directions $\mathbf{a}^{1}$ and $\mathbf{a}^{2}$.}

\bigskip

This theorem implies that for solving the interpolation problem on $L$ it is
enough to verify whether $L$ contains a closed path. We use this fact in the
proof of Theorem 1.1. In other words, we show in the next section that a set
of distinct three lines in $\mathbb{R}^{n}$ possesses closed paths with
respect to two arbitrarily given directions.

\

\begin{center}
{\large \textbf{2. The three-line interpolation problem}}
\end{center}

The following theorem is equivalent to Theorem 1.1.

\bigskip

\textbf{Theorem 2.1.} \textit{Assume we are given two linearly independent
directions $\mathbf{a}^{1},\mathbf{a}^{2}\in \mathbb{R}^{n}$. Assume, in
addition, $l_{i}=\{t\mathbf{b}^{i}+\mathbf{c}^{i}:t\in \mathbb{R}\},$ $%
i=1,2,3,$ are distinct straight lines in $\mathbb{R}^{n}$ and $L=l_{1}\cup
l_{2}\cup l_{3}$. Then there exists a closed path $p=\left( \mathbf{x}%
^{1},...,\mathbf{x}^{m}\right) $ with respect to the directions $\mathbf{a}%
^{1}$ and $\mathbf{a}^{2}$ such that $p\subset L$.}

\bigskip

\textbf{Proof.} Complete the pair $\{\mathbf{a}^{1},\mathbf{a}^{2}\}$ to a
basis $\{\mathbf{a}^{1},\mathbf{a}^{2},...,\mathbf{a}^{n}\}$ in $\mathbb{R}%
^{n}$. The nonsingular linear transformation of coordinates $S:\mathbf{x}%
=(x_{1},...,x_{n})\rightarrow \mathbf{y}=(y_{1},...,y_{n}),$ $y_{i}=\mathbf{a%
}^{i}\cdot $ $\mathbf{x},$ $i=1,...,n$, has the property $\mathbf{e}%
^{i}\cdot $ $\mathbf{y}=\mathbf{a}^{i}\cdot $ $\mathbf{x}$, where $\mathbf{e}%
^{i}$ is the $i$-th unit vector (that is, the vector having $1$ in the $i$%
-th position and zeros everywhere else). Therefore, if $\left( \mathbf{x}%
^{1},...,\mathbf{x}^{m}\right) $ is a path with respect to $\mathbf{a}^{1},%
\mathbf{a}^{2}$, then the transformed set $\left( \mathbf{y}^{1},...,\mathbf{%
y}^{m}\right) $ will be a path with respect to $\mathbf{e}^{1},\mathbf{e}%
^{2} $. This fact admits to choose $\mathbf{a}^{1}=\mathbf{e}%
^{1}=(1,0,...,0) $ and $\mathbf{a}^{2}=\mathbf{e}^{2}=(0,1,...,0)$. Thus we
will prove the theorem if we prove it for the basis directions $\mathbf{e}%
^{1}$ and $\mathbf{e}^{2}$.

Let $P$ be the orthogonal projection of $\mathbb{R}^{n}$ onto the plane
generated by the vectors $\mathbf{e}^{1}$ and $\mathbf{e}^{2}$, that is,
onto the $x_{1}x_{2}$ plane. Clearly, the image of a straight line under
this projection will be a straight line or a point. Note that if for some $%
l_{i}$ the image of it is a point, then any two distinct points $\mathbf{x}%
^{1},\mathbf{x}^{2}\in l_{i}$ form a closed path, since in this case $%
\mathbf{x}^{2}-\mathbf{x}^{1}$ is orthogonal to both $\mathbf{e}^{1}$ and $%
\mathbf{e}^{2}$. Further note that the projections of two (or all three) of $%
l_{i}$ may coincide. In this case we also have a closed path $(\mathbf{x}%
^{1},\mathbf{x}^{2})$ with the vertices $\mathbf{x}^{1}$ and $\mathbf{x}^{2}$
lying in different straight lines. We see that it is enough to consider the
general case when the images of all the lines $l_{i}$ are again straight
lines and these images do not coincide. If we prove that the image $%
PL=Pl_{1}\cup Pl_{2}\cup Pl_{3}\subset \mathbb{R}^{2}$ contains a closed
path $\left( \mathbf{x}^{1},...,\mathbf{x}^{2m}\right) $ with respect to the
vectors $(1,0)$ and $(0,1)$, then the initial set $L$ contains a closed path
with respect to $\mathbf{e}^{1}$ and $\mathbf{e}^{2}$. Indeed, if for some $%
\mathbf{x}^{i}$, $i=1,...,2m$, the preimage $P^{-1}\mathbf{x}^{i}\cap L$ has
more than a single point, then any two of them form a closed path. And if $%
\left( \mathbf{x}^{1},...,\mathbf{x}^{2m}\right) $ is a closed path in $PL$
and for each $\mathbf{x}^{i}$, the preimage $P^{-1}\mathbf{x}^{i}\cap L$ is
a single point, denoted here by $\mathbf{z}^{i}$, then $\left( \mathbf{z}%
^{1},...,\mathbf{z}^{2m}\right) $ is a closed path in $L.$ This observation
reduces the problem to the two-dimensional case and we will prove the
theorem if we prove that the union of any distinct three lines in $\mathbb{R}%
^{2}$ contain a closed path with respect to the coordinate directions.

\bigskip

Concerning three lines in $\mathbb{R}^{2}$ we have only the following
options.

\bigskip

\textbf{Case 1.} At least one of $l_{i}$ is parallel to the $x_{1}$ or $%
x_{2} $ axis. In this case it is not difficult to see that there exist
four-point closed paths.

\bigskip

\textbf{Case 2.} The straight lines $l_{i}=\{t\mathbf{b}^{i}+\mathbf{c}%
^{i}:t\in \mathbb{R}\},$ $i=1,2,3,$ are parallel. That is, $\mathbf{b}^{1}=%
\mathbf{b}^{2}=\mathbf{b}^{3}.$ Let $\mathbf{b}^{i}=\mathbf{b}=(b_{1},b_{2})$%
, $\mathbf{c}^{i}=(c_{1}^{i},c_{2}^{i}),$ $i=1,2,3$.

\textit{Linear algebraic approach. }Consider a six-point set $(\mathbf{x}%
^{1},\mathbf{x}^{2},\mathbf{x}^{3},\mathbf{x}^{4},\mathbf{x}^{5},\mathbf{x}%
^{6})$ in $l_{1}\cup l_{2}\cup l_{3}$ with the property that $\mathbf{x}^{1},%
\mathbf{x}^{4}\in l_{1},$ $\mathbf{x}^{2},\mathbf{x}^{5}\in l_{2},$ and $%
\mathbf{x}^{3},\mathbf{x}^{6}\in l_{3}$. Put $\mathbf{x}^{k}=t_{k}\mathbf{b}+%
\mathbf{c}^{k},$ $k=1,...,6$, $\mathbf{c}^{j+3}=\mathbf{c}^{j}$, $j=1,2,3.$
Let us check if these points could form a closed path in the given order.
This is possible if the system of linear equations

\begin{equation*}
\left\{
\begin{array}{c}
t_{1}b_{1}+c_{1}^{1}=t_{2}b_{1}+c_{1}^{2} \\
t_{2}b_{2}+c_{2}^{2}=t_{3}b_{2}+c_{2}^{3} \\
t_{3}b_{1}+c_{1}^{3}=t_{4}b_{1}+c_{1}^{1} \\
t_{4}b_{2}+c_{2}^{1}=t_{5}b_{2}+c_{2}^{2} \\
t_{5}b_{1}+c_{1}^{2}=t_{6}b_{1}+c_{1}^{3} \\
t_{6}b_{2}+c_{2}^{3}=t_{1}b_{2}+c_{2}^{1}%
\end{array}%
\right. \eqno(2.1)
\end{equation*}%
has a solution. In fact, applying Cramer's rule, one can see that this
system has infinitely many solutions $(t_{1},t_{2},...,t_{6})$. This means
that we have infinitely many six-point closed paths with vertices positioned
at these lines. Note that if $(t_{1},t_{2},...,t_{6})$ is a solution of the
system (2.1), then $(t_{1}+\alpha ,t_{2}+\alpha ,...,t_{6}+\alpha )$ is also
a solution for any $\alpha \in \mathbb{R}$. Thus, we see that all six-point
closed paths satisfying (2.1) can be obtained by a single six-point closed
path by the translation $T(\mathbf{x})=\mathbf{x}+\alpha \mathbf{b}$.

\textit{Geometric approach. }There is also a geometric approach to this
case. Consider a path $(A,B,C,D,E,F,A^{\prime })$ such that $A,D,A^{\prime
}\in l_{1}$, $B,E\in l_{2}$, $C,F\in l_{3}$. That is, the vertices of this
path lie alternatively at the lines $l_{1},l_{2},l_{3}$. It is not difficult
to see that $\overrightarrow{AB}+\overrightarrow{CD}+\overrightarrow{EF}=%
\overrightarrow{\mathbf{0}}$ and $\overrightarrow{BC}+\overrightarrow{DE}+%
\overrightarrow{FA^{\prime }}=\overrightarrow{\mathbf{0}}$. Hence $%
\overrightarrow{AA^{\prime }}=\overrightarrow{\mathbf{0}}$ which yields $%
A=A^{\prime }$. Thus the path $(A,B,C,D,E,F)$ is closed. See Fig. 2 for the
geometric illustration of this case.

\begin{figure}[ht]
\hfill
\par
\begin{center}
\includegraphics[width=0.65\textwidth]{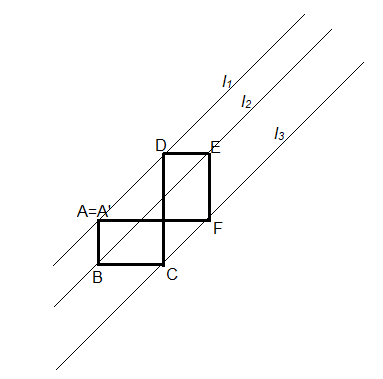}
\end{center}
\caption{\textit{The six-point closed path $(A,B,C,D,E,F)$ in $l_{1}\cup
l_{2}\cup l_{3}$.}}
\end{figure}

\bigskip

\textbf{Case 3.} The lines $l_{i}$ have a common point of intersection $(\xi
_{1},\xi _{2})$ and none of them is parallel to $x_{1}$ or $x_{2}$ axis. By
translation we can make that $(\xi _{1},\xi _{2})=(0,0)$. So without loss of
generality we may assume that the lines intersect at the coordinate origin.
Let $x_{2}=px_{1}$, $x_{2}=qx_{1}$, $x_{2}=rx_{1}$ be equations of these
three lines. Then it is easy to verify that the set

\begin{equation*}
\{(p,pq),(p,pr),(r,pr),(r,qr),(q,qr),(q,pq)\}
\end{equation*}%
is a closed path in $l_{1}\cup l_{2}\cup l_{3}$.

\bigskip

\textbf{Case 4.} The lines do not intersect at the same point, they are not
parallel and none of them is parallel to the coordinate axis. We claim that
in this case there is a four-point closed path $\left\{
(u_{1},u_{2}),(u_{1},v_{2}),(v_{1},v_{2}),(v_{1},u_{2})\right\} $ in the
union $L$ of these lines. Let us prove this. Assume $x_{2}=ax_{1}+b,$ $%
x_{2}=cx_{1}+d,$ $x_{2}=ex_{1}+f$ are the equations of $l_{1},l_{2},l_{3}$,
respectively.

\textit{Linear algebraic approach.} We seek such a closed path with the
vertices $(u_{1},u_{2}),(v_{1},v_{2})$ in one line and with $%
(u_{1},v_{2}),(v_{1},u_{2})$ each in other two lines, respectively. If $%
(u_{1},u_{2}),(v_{1},v_{2})\in l_{1}$, $(u_{1},v_{2})\in l_{2}$, $%
(v_{1},u_{2})\in l_{3}$, then we have the system of linear equations

\begin{equation*}
\left\{
\begin{array}{c}
-au_{1}+u_{2}=b \\
-av_{1}+v_{2}=b \\
-cu_{1}+v_{2}=d \\
-ev_{1}+u_{2}=f%
\end{array}%
\right. .
\end{equation*}%
The determinant of this system (up to a sign) is equal to $ce-a^{2}$. In
other two cases when $(u_{1},u_{2}),(v_{1},v_{2})\in l_{2}$, $%
(u_{1},v_{2})\in l_{1}$, $(v_{1},u_{2})\in l_{3}$ and when $%
(u_{1},u_{2}),(v_{1},v_{2})\in l_{3}$, $(u_{1},v_{2})\in l_{1}$, $%
(v_{1},u_{2})\in l_{2}$, the determinants of the corresponding systems will
be $ae-c^{2}$ and $ac-e^{2}$, respectively. If the expressions $ce-a^{2}$, $%
ae-c^{2}$ and $ac-e^{2}$ are all zero, then $a=c=e$ and hence the lines are
parallel. This reduces the problem to Case 2. If one of these expressions is
not zero, then the corresponding system has a unique solution $%
u_{1},u_{2},v_{1},v_{2}$. Note that $u_{1}=v_{1}$ or $u_{2}=v_{2}$ is
possible if and only if the lines intersect at the same point, which reduces
the problem to Case 3 above. Thus we obtain that there exists a four-point
closed path in the union of the given three straight lines. Certainly, the
vertices of this path can be explicitly determined by the coefficients $%
a,b,c,d,e,f$.

\textit{Geometric approach.} Let $\Omega $ be the set of all non-degenerate
rectangles with sides parallel to the coordinate axis. Let $l_{1}$ and $%
l_{2} $ intersect at a point $O$ forming four rays $r_{1},r_{2},r_{3},r_{4}$
with a vertex at $O$. If $r$ is any of these rays and $l_{i}$ is the line
not containing $r$, then for any $P\in r\backslash \{O\}$, we let $\Pi
_{P}\in \Omega $ be the rectangle with one vertex at $P$ and two adjacent
(to $P$) vertices at $l_{i}$. We will denote the fourth vertex of this
rectangle by $v(P)$. Then the set $\{v(P):P\in r\backslash \{O\}\}\cup \{O\}$
is a ray at the vertex $O$; If $r=r_{i}$, then we denote this ray by $s_{i}$. The rays $s_{1}, s_{2}, s_{3}, s_{4}$ are pairwise distinct and lie on two different lines passing through $O$. Since $l_{3}$ cannot be parallel to both of these lines, it must intersect at least one of the rays $s_{i}$. This intersection point $v(P)\in l_{3}$ and the
three points considered above form a four-closed path in $l_{1}\cup
l_{2}\cup l_{3} $. See Fig. 3 for the geometric illustration of this case.

\begin{figure}[ht]
\hfill
\par
\begin{center}
\includegraphics[width=0.75\textwidth]{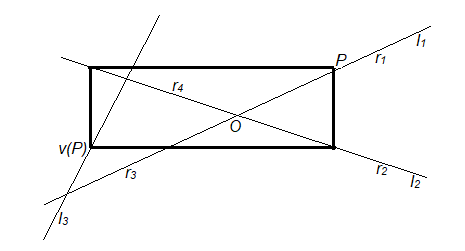}
\end{center}
\caption{\textit{Four-point closed path in $l_{1}\cup l_{2}\cup l_{3}$. The
vertex $v(P)$ depends on $P$. There exists $P \in l_{1}\cup l_{2}$ such that
$v(P) \in l_{3}$.}}
\end{figure}

Theorem 2.1 has been proven.

\bigskip

\textbf{Remark.} One and two straight lines in $\mathbb{R}^{n}$ may not
contain closed paths, hence making the ridge function interpolation
possible. The characterization of such straight lines was given in \cite{IP}.

\

\begin{center}
{\large \textbf{3. The case of three or more directions}}
\end{center}

Let us consider the interpolation problem in the case when we have three or
more directions. Given $m$ straight lines $l_{1},...,l_{m}$ in $\mathbb{R}%
^{n}$, can we interpolate arbitrary data on these lines by functions from $%
\mathcal{R}\left( \mathbf{a}^{1},...,\mathbf{a}^{r}\right) $, $r\geq 3$? In
\cite{IP}, it was conjectured that it should be possible, except in certain
specific cases, to interpolate along $m$ straight lines, if $m\leq r$. And
it should be impossible to interpolate arbitrary data on any $m$ straight
lines, if $m>r$. Although, we do not know a solution to this general problem
of interpolation, we want as in the previous sections to reduce it to the
more tractable problem. The latter will involve cycles -- mathematical
objects which are generalizations of closed paths (see \cite[Chapter 1]{Ism}%
). Let $\delta _{A}$ denote the characteristic function of a set $\ A\subset
\mathbb{R}.$ That is,

\begin{equation*}
\delta _{A}(y)=\left\{
\begin{array}{c}
1,~if~y\in A \\
0,~if~y\notin A.%
\end{array}%
\right.
\end{equation*}%
A set of points $\{\mathbf{x}^{1},\ldots ,\mathbf{x}^{k}\}\subset \mathbb{R}%
^{n}$ is called a \textit{cycle} with respect to the directions $\mathbf{a}%
^{1},...,\mathbf{a}^{r}$ if there exists a vector $\lambda =(\lambda
_{1},\ldots ,\lambda _{k})$ with the nonzero components such that for any $%
i=1,\ldots ,r$%
\begin{equation*}
\sum_{j=1}^{k}\lambda _{j}\delta _{\mathbf{a}^{i}\cdot \mathbf{x}^{j}}(t)=0,%
\eqno(3.1)
\end{equation*}
for all $t\in \mathbb{R}$.

Let us explain Eq. (3.1) in detail. We will see that, in fact, it stands for
a system of linear equations. Fix the subscript $i.$ Let the set $\{\mathbf{a%
}^{i}\cdot \mathbf{x}^{j},$ $j=1,...,k\}$ have $s_{i}$ different values,
which we denote by $\gamma _{1}^{i},\gamma _{2}^{i},...,\gamma _{s_{i}}^{i}.$
Take the first number $\gamma _{1}^{i}.$ Putting $t=\gamma _{1}^{i}$, we
obtain from (3.1) that

\begin{equation*}
\sum_{j}\lambda _{j}=0,
\end{equation*}%
where the sum is taken over all $j$ such that $\mathbf{a}^{i}\cdot \mathbf{x}%
^{j}=\gamma _{1}^{i}.$ This is the first linear equation in $\lambda _{j}$
corresponding to $\gamma _{1}^{i}$. Take now $\gamma _{2}^{i}$. By the same
way, putting $t=\gamma _{2}^{i}$ in (3.1), we can form the second equation.
Continuing until $\gamma _{s_{i}}^{i}$, we obtain $s_{i}$ linear homogeneous
equations in $\lambda _{1},...,\lambda _{k}$. The coefficients of these
equations are the integers $0$ and $1$. By varying $i$, we finally obtain $%
s=\sum_{i=1}^{r}s_{i}$ such equations. Hence (3.1), in its expanded form,
stands for the system of these linear equations. Thus $\{\mathbf{x}%
^{1},\ldots ,\mathbf{x}^{k}\}$ is a cycle if the system of linear equations
of the form (3.1) has a solution with nonzero components.

For example, assume $r=2,$ $\mathbf{a}^{1}\cdot \mathbf{x}^{1}=\mathbf{a}%
^{1}\cdot \mathbf{x}^{2}$, $\mathbf{a}^{2}\cdot \mathbf{x}^{2}=\mathbf{a}%
^{2}\cdot \mathbf{x}^{3}$, $\mathbf{a}^{1}\cdot \mathbf{x}^{3}=\mathbf{a}%
^{1}\cdot \mathbf{x}^{4}$,..., $\mathbf{a}^{2}\cdot \mathbf{x}^{k-1}=\mathbf{%
a}^{2}\cdot \mathbf{x}^{k}$,$\mathbf{a}^{1}\cdot \mathbf{x}^{k}=\mathbf{a}%
^{1}\cdot \mathbf{x}^{1}$ Then it is not difficult to see that for a vector $%
\lambda =(\lambda _{1},\ldots ,\lambda _{k})$ with the components $\lambda
_{i}=(-1)^{i},$ we have
\begin{equation*}
\sum_{j=1}^{k}\lambda _{j}\delta _{\mathbf{a}^{i}\cdot \mathbf{x}^{j}}=0,
\end{equation*}%
for $i=1,2.$ Thus, the set $p=\{\mathbf{x}^{1},\ldots ,\mathbf{x}^{k}\}$ is
a cycle with respect to the directions $\mathbf{a}^{1}$ and $\mathbf{a}^{2}$%
. Note that $p$ is also a closed path (after some suitable permutation of
its points which we assume to be as given). Hence any closed path is a
cycle. It is not difficult to see that any cycle with respect to two
directions is a union of closed paths with respect to the same directions.

Another example is the set $\{(0,0,0),~(0,0,1),~(0,1,0),~(1,0,0),~(1,1,1)\}$
in $\mathbb{R}^{3}$, which forms a cycle with respect to the coordinate
directions. In this example, the vector $\lambda $ above can be taken as $%
(-2,1,1,1,-1).$

In \cite{Is1}, the second author, in particular, proved that any function $F$
defined on a given $X\subset \mathbb{R}^{n}$ belongs to $\mathcal{R}\left(
\mathbf{a}^{1},...,\mathbf{a}^{r}\right) $ if and only if $X$ does not
contain cycles with respect to $\mathbf{a}^{1},...,\mathbf{a}^{r}$. Thus the
interpolation problem on $m$ straight lines can be reduced to the problem of
existence of cycles in these lines. In fact, the interpolation problem is
solvable if and only if the union of these lines does not contain cycles.
Paralleling Theorem 2.1 above, we speculate that any given $m$ distinct
straight lines in $\mathbb{R}^{n}$ contain a cycle with respect to $r$
directions if and only if $r<m$. This is an equivalent formulation of the
above conjecture.

\bigskip

\textbf{Acknowledgments.} The second author is grateful to Namig Guliyev for
inspiring discussions and his input to the analysis of Case 4 above.

\end{document}